\documentstyle[12pt]{article}
\newcommand{\const}{\mathop{\rm const}\limits}

\newcommand{\supp}{\mathop{\rm supp}\limits}

\newcommand{\card}{\mathop{\rm card}\limits}

\newcommand{\Var}{\mathop{\rm Var}\limits}

\newcommand{\sign}{\mathop{\rm sign}\limits}

\newcommand{\Law}{\mathop{\rm Law}\limits}

\newcommand{\Sub}{\mathop{\rm Sub}\limits}

\newcommand{\SSub}{\mathop{\rm SSub}\limits}

\newcommand{\grad}{\mathop{\rm grad}\limits}

\begin{document}

\begin{center}

{\bf  VECTOR REARRANGEMENT INVARIANT BANACH SPACES \\

\vspace{3mm}

 OF RANDOM VARIABLES WITH EXPONENTIAL \\

\vspace{3mm}

DECREASING TAILS OF DISTRIBUTIONS }\\

\vspace{5mm}

 {\sc Ostrovsky E., Sirota L.}\\

\vspace{3mm}

 \ Department of Mathematics and Statistics, Bar-Ilan University, \\
59200, Ramat Gan, Israel.\\
e-mail: \ eugostrovsky@list.ru \\

\vspace{3mm}

 \ Department of Mathematics and Statistics, Bar-Ilan University,\\
59200, Ramat Gan, Israel.\\
e-mail: \ sirota3@bezeqint.net \\

\vspace{3mm}

 {\sc Abstract.}

\end{center}

 \  We present in this paper the theory of multivariate Banach spaces of random variables
with exponential decreasing tails of distributions. \par

 \vspace{3mm}

 {\it Key words and phrases:}  Random variable and random vector (r.v.), centered (mean zero) r.v., moment generating function,
rearrangement invariant Banach space on vector random variables (vectors), ordinary and exponential moments, absolutely and
relatively monotonic functions, multivariate Bernstein's theorem, Grand Lebesgue Space (GLS) norm, octant,
tail of distribution, Young-Orlicz function, norm, Chernov's estimate, theorem of Fenchel-Moraux, generating function, upper and lower
estimates, non-asymptotical exponential estimations, Kramer's condition.

\vspace{4mm}

{\it Mathematics Subject Classification (2000):} primary 60G17; \ secondary
 60E07; 60G70.\\

\vspace{5mm}

\section{ Introduction. Previous results.} \par

\vspace{4mm}

 \ We present here for beginning some known facts from the theory of one-dimensional random variables
with exponential decreasing tails of distributions, see  \cite{Ostrovsky4},  \cite{Kozachenko1}, \cite{Ostrovsky1}, chapters 1,2. \par

 \ Especially  we menton the authors preprint \cite{Ostrovsky5}; we offer in comparison with existing there results a more fine approach. \par

\vspace{3mm}

 \ Let $ (\Omega,F,{\bf P} ) $ be a probability space, $ \Omega = \{\omega\}. $ \par

\vspace{3mm}

 \ Let also $ \phi = \phi(\lambda), \lambda \in (-\lambda_0, \lambda_0), \ \lambda_0 =
\const \in (0, \infty] $ be certain even strong convex which takes positive values for positive arguments twice continuous
differentiable function, such that
$$
 \phi(0) = \phi'(0) = 0, \ \phi^{''}(0) > 0, \ \lim_{\lambda \to \lambda_0} \phi(\lambda)/\lambda = \infty. \eqno(1.1)
$$
  \ We denote the set of all these (Young-Orlicz) function as $ \Phi; \ \Phi = \{ \phi(\cdot) \}. $ \par

 \ We say by definition that the {\it centered} random variable (r.v) $ \xi = \xi(\omega) $
belongs to the space $ B(\phi), $ if there exists some non-negative constant
$ \tau \ge 0 $ such that

$$
\forall \lambda \in (-\lambda_0, \lambda_0) \ \Rightarrow
\max_{\pm} {\bf E} \exp(\pm \lambda \xi) \le \exp[ \phi(\lambda \ \tau) ]. \eqno(1.2)
$$

 \ The minimal value $ \tau $ satisfying (1.2) for all the values $  \lambda \in (-\lambda_0, \lambda_0), $
is named a $ B(\phi) \ $ norm of the variable $ \xi, $ write

$$
||\xi||B(\phi)  \stackrel{def}{=}
$$

 $$
 \inf \{ \tau, \ \tau > 0: \ \forall \lambda:  \ |\lambda| < \lambda_0 \ \Rightarrow
  \max_{\pm}{\bf E}\exp( \pm \lambda \xi) \le \exp(\phi(\lambda \ \tau)) \}. \eqno(1.3)
 $$

 \ These spaces are very convenient for the investigation of the r.v. having a
exponential decreasing tail of distribution, for instance, for investigation of the limit theorem,
the exponential bounds of distribution for sums of random variables,
non-asymptotical properties, problem of continuous and weak compactness of random fields,
study of Central Limit Theorem in the Banach space etc.\par

  The space $ B(\phi) $ with respect to the norm $ || \cdot ||B(\phi) $ and
ordinary algebraic operations is a Banach space which is isomorphic to the subspace
consisted on all the centered variables of Orlicz's space $ (\Omega,F,{\bf P}), N(\cdot) $
with $ N \ - $ function

$$
N(u) = \exp(\phi^*(u)) - 1, \ \phi^*(u) = \sup_{\lambda} (\lambda u -
\phi(\lambda)).
$$
 \ The transform $ \phi \to \phi^* $ is called Young-Fenchel transform. The proof of considered
assertion used the properties of saddle-point method and theorem of Fenchel-Moraux:
$$
\phi^{**} = \phi.
$$

  \ The next facts about the $ B(\phi) $ spaces are proved in \cite{Kozachenko1}, \cite{Ostrovsky1}, p. 19-40:

$$
{\bf 1.} \ \xi \in B(\phi) \Leftrightarrow {\bf E } \xi = 0, \ {\bf and} \ \exists C = \const > 0,
$$

$$
U(\xi,x) \le \exp(-\phi^*(Cx)), x \ge 0, \eqno(1.4)
$$
where $ U(\xi,x)$ denotes in this section the {\it one-dimensional tail} of
distribution of the r.v. $ \xi: $

$$
U(\xi,x) = \max \left( {\bf P}(\xi > x), \ {\bf P}(\xi < - x) \right), \ x \ge 0,
$$
and this estimation is in general case asymptotically as $ x \to \infty  $ exact. \par

 \ Here and further $ C, C_j, C(i) $ will denote the non-essentially positive
finite "constructive" constants.\par
 More exactly, if $ \lambda_0 = \infty, $ then the following implication holds:

$$
\lim_{\lambda \to \infty} \phi^{-1}(\log {\bf E} \exp(\lambda \xi))/\lambda =
K \in (0, \infty)
$$
if and only if

$$
\lim_{x \to \infty} (\phi^*)^{-1}( |\log U(\xi,x)| )/x = 1/K.
$$

 \ Hereafter $ f^{-1}(\cdot) $ denotes the inverse function to the
function $ f $ on the left-side half-line $ (C, \infty). $ \par

  \ Let $  F =  \{  \xi(t) \}, \ t \in T, \ T  $ is an arbitrary set, be the family of somehow
dependent mean zero random variables. The function $  \phi(\cdot) $ may be "constructive" introduced by the formula

$$
\phi(\lambda) = \phi_F(\lambda) \stackrel{def}{=} \max_{\pm} \log \sup_{t \in T}
 {\bf E} \exp(  \pm \lambda \xi(t)), \eqno(1.5)
$$
 if obviously the family $  F  $ of the centered r.v. $ \{ \xi(t), \ t \in T \} $ satisfies the  so-called
{\it uniform } Kramer's condition:
$$
\exists \mu \in (0, \infty), \ \sup_{t \in T} U(\xi(t), \ x) \le \exp(-\mu \ x),
\ x \ge 0.
$$
 In this case, i.e. in the case the choice the function $ \phi(\cdot) $ by the
formula (1.5), we will call the function $ \phi(\lambda) = \phi_0(\lambda) $
a {\it natural } function. \par

\vspace{3mm}

 \ {\bf 2.} We define $ \psi(p) = p/\phi^{-1}(p), \ p \ge 2. $ \par

\vspace{3mm}

 \ Let us introduce a new norm, the so-called "moment norm", or equally Grand Lebesgue Space (GLS) norm,
on the set of r.v. defined in our probability space by the following way: the
space $ G(\psi) $ consist, by definition, on all the centered (mean zero) r.v. with finite norm

$$
||\xi||G(\psi) \stackrel{def}{=} \sup_{p \ge 1} \left[ |\xi|_p/\psi(p) \right],  \eqno(1.6)
$$
here and in what follows as ordinary

$$
\ |\xi|_p :={\bf E}^{1/p} |\xi|^p = \left[  \int_{\Omega} |\xi(\omega)|^p \ {\bf P}(d \omega)  \right]^{1/p}.
$$

 \ It is proved that the spaces $ B(\phi) $ and $ G(\psi) $ coincides: $ B(\phi) =
G(\psi) $ (set equality) and both
the norms $ ||\cdot||B(\phi) $ and $ ||\cdot|| $ are linear equivalent: $ \exists C_1 =
C_1(\phi), C_2 = C_2(\phi) = \const \in (0,\infty), \ \forall \xi \in B(\phi) \ \Rightarrow $

$$
||\xi||G(\psi) \le C_1 \ ||\xi||B(\phi) \le C_2 \ ||\xi||G(\psi). \eqno(1.6a)
$$

\vspace{3mm}

 \ {\bf 3.} The definition (1.6) is correct still for the non-centered random
variables $ \xi.$ If for some non-zero r.v. $ \xi \ $ we have $ ||\xi||G(\psi) < \infty, $ then for all positive values $ u $

$$
{\bf P}(|\xi| > u) \le 2 \ \exp \left( - u/(C_3 \ ||\xi||G(\psi)) \right).
\eqno(1.7)
$$
and conversely if a r.v. $ \xi $ satisfies (1.7), then $ ||\xi||G(\psi) < \infty. $ \par

\vspace{4mm}

 {\bf  We intend to extend in this report the definition and properties of these spaces into a
multidimensional case, i.e. when all the considered r.v. are vectors. } \par

\vspace{4mm}

 \ We agree to omit the proofs for the multivariate case if they are likewise to ones  for the ordinary
numerical variables.

\vspace{4mm}

 \ The paper is organized as follows. In the second section we give necessary definitions, conventions and restrictions.
 We prove in the section 3 the Banach structure of offered spaces. The fourth section is devoted to the description of the
characterization features for multidimensional moment generating function. \par

 \ In the fifth section we investigate the linear transforms for vectors in these spaces.
 In the next section we obtain some multidimensional bilateral tail inequalities for random
vectors belonging to these spaces. The seventh section contains states the
equivalence with moments (Grand Lebesgue Spaces) and introduced norms.\par

 \ We consider in the $   8^{th} $ one of the important applications of introduced spaces: the
exponential bounds for the sums of independent random vectors. We discover in the penultimate section the
relation between our spaces with multivariate Orlicz spaces. \par

 \ The last section contains some concluding remarks.\par

\vspace{4mm}

\section{ Notations. Definitions, conventions and restrictions. } \par

\vspace{4mm}

 \ We need to introduce some new notations. \par

\vspace{3mm}

 \ Denote by $  \epsilon  = \vec{\epsilon} = \{ \epsilon(1),   \epsilon(2), \ldots,   \epsilon(d) \} $ the non-random
$  d \ - $ dimensional  numerical vector, $  d = 2,3,\ldots, $ whose components take the values $  \pm 1 $ only. \par
 \ Set $  \vec{1} = (1,1,\ldots,1) \in R^d_+. $ \par

 \ Denote by $  \Theta = \Theta(d) = \{ \ \vec{\epsilon} \ \} $ {\it collection} of all such a vectors.  Note that
 $  \card \Theta = 2^d $  and $ \vec{1} \in \Theta.$  \par

 \ Another denotations. For $ \vec{\epsilon} \in \Theta(d) $  and vector $ \vec{x} $ we introduce the coordinatewise product
as a $ d  \ - $ dimensional vector of the form

$$
\vec{\epsilon} \otimes \vec{x} \stackrel{def}{=} \{ \epsilon(1) \ x(1), \ \epsilon(2) \ x(2), \ \ldots,  \epsilon(d) \ x(d)  \}.
$$

 \ Let $ f = f(\lambda), \ \lambda \in R^d $ be sufficiently smooths function  and
  $ \vec{k} $ be non-negative $ d  \ - $ dimensional integer  vectors $ \vec{k} = \{ k(1), \ k(2), \ldots, k(d)  \};
  |k| := \sum_j k(j). $  We denote as ordinary

$$
\frac{\partial^{ \vec{k}} \ f }{ \partial \vec{\lambda}^{{\vec k}}} =
\frac{\partial^{|k|} \ f(\lambda)}{ \partial \lambda(1)^{k(1)} \ \partial \lambda(2)^{k(2)} \  \ldots \partial \lambda(d)^{k(d) }}.
\eqno(2.0)
$$

\vspace{3mm}

 \ {\bf Definition 2.1.}\par

\vspace{3mm}

 \ Let $ \xi = \vec{\xi} = (\xi(1), \xi(2), \ldots, \xi(d) )  $ be a centered random vector such that each its component
$  \xi(j) $  satisfies the Kramer's condition. The {\it natural function} $ \phi_{\xi}= \phi_{\xi}(\lambda), \
 \lambda = \vec{\lambda} = (\lambda(1), \lambda(2), \ldots, \lambda(d))  \in R^d  $
for the random vector $  \xi  $ is defined as follows:

$$
\exp \{\phi_{\xi}(\lambda)\} \stackrel{def}{=} \max_{\vec{\epsilon}} \
{\bf E} \exp \left\{ \sum_{j=1}^d  \epsilon(j) \lambda(j) \xi(j) \right\}=
$$

$$
\max_{\vec{\epsilon} \in \Theta} \
{\bf E} \exp \{\epsilon(1) \lambda(1) \xi(1) + \epsilon(2) \lambda(2)\xi(2)+ \ldots +   \epsilon(d) \lambda(d)\xi(d) \}, \eqno(2.1)
$$
 where $  "\max" $ is calculated over all the combinations of signs $ \epsilon(j) =  \pm 1. $\par

\vspace{3mm}

 \ {\bf Definition 2.2.}\par

\vspace{3mm}

 \ The {\it tail function } for the random vector $ \vec{\xi} \ $
   $ U(\vec{\xi}, \vec{x}), \ \vec{x} = (x(1), x(2), \ldots, x(d)),   $ where all the coordinates $ x(j)  $
  of the deterministic vector $ \vec{x} $ are non-negative, is defined as follows.

$$
U(\vec{\xi}, \vec{x}) \stackrel{def}{=} \max_{\vec{\epsilon}}
{\bf P} \left( \cap_{j=1}^d  \{  \epsilon(j) \xi(j) > x(j) \}  \right) =
$$

$$
 \max_{ \vec{\epsilon} \in \Theta}
{\bf P}(\epsilon(1) \xi(1) > x(1), \ \epsilon(2) \xi(2) > x(2), \ \ldots, \ \epsilon(d) \xi(d) > x(d) ), \eqno(2.2)
$$
 where as before $  "\max" $ is calculated over all the combinations of signs $ \epsilon(j) =  \pm 1. $\par

 \ We illustrate this notion in the case $  d = 2.  $ Let $ \vec{\xi} = (\xi(1), \ \xi(2))  $ be a two-dimensional random
vector and let $ x, y  $  be non-negative numbers. Then

$$
U( (\xi(1), \xi(2)), \ (x,y) ) =
$$

$$
\max [ {\bf P} (\xi(1) > x, \ \xi(2) > y), \ {\bf P} (\xi(1) > x, \ \xi(2) < - y),
$$

$$
{\bf P} (\xi(1) < - x, \ \xi(2) > y), \ {\bf P} (\xi(1) < - x, \ \xi(2) < - y) ].
$$

\vspace{3mm}

 \ {\bf Definition 2.3.}\par

\vspace{3mm}

 \ Let $  h = h(x), \ x \in R^d  $ be some non-negative real valued function, which is finite on some non-empty
neighborhood of origin.  We denote as ordinary

$$
\supp h = \{x, \ h(x) < \infty   \}.
$$

 \ The Young-Fenchel, or Legendre transform $  h^*(y), \ y \in R^d  $ is defined likewise the one-dimensional case

$$
h^*(y) \stackrel{def}{=} \sup_{x \in \supp h} ( (x,y) - h(x)). \eqno(2.3)
$$

 \ Herewith $ (x,y)  $ denotes the scalar product of the vectors $ x,y: \ (x,y) = \sum_j x(j)y(j);
  \ |x| = \sqrt{(x,x)}. $\par

 \ Obviously, if the set $ \supp h $ is central symmetric, then the function $ h^*(y) $ is even.\par

\vspace{3mm}

 \ {\bf Definition 2.4.}\par

\vspace{3mm}

 \ Recall, see \cite{Rao1}, \cite{Rao2} that the function $  x \to g(x), \ x \in R^d, \ g(x) \in R^1_+  $ is named
 multivariate Young, or Young-Orlicz function, if it is even,  $ d \ - $ times  continuous differentiable, convex, non-negative,
 finite on the whole space $  R^d,$ and such that

$$
g(x) = 0 \ \Leftrightarrow x = 0; \hspace{4mm} \frac{\partial g}{\partial x}/(\vec{x} = 0) = 0,
$$

$$
\det \frac{\partial^2 g}{\partial x^2} /(\vec{x} = 0) > 0. \eqno(2.4)
$$

 We explain in detail:

$$
 \frac{\partial g}{\partial x} = \left\{  \frac{ \partial g}{ \partial x_j} \right\} = \grad g,  \hspace{5mm}
\frac{\partial^2 g}{\partial x^2} = \left \{ \frac{\partial^2 g}{\partial x_k \partial x_l}  \right\}.
\eqno(2.5)
$$

 \ We assume  finally

$$
\lim_{|x| \to \infty} \frac{\partial^d g }{\prod_{k=1}^d  \partial x_k} = \infty. \eqno(2.6)
$$

 \ We will denote the set of all such a functions by $  Y = Y(R^d) $ and denote also by $ D  $ introduced before matrix

$$
D = D_g := \frac{1}{2} \left\{ \frac{\partial^2 g(0)}{\partial x_k \partial x_l} \right\}.
$$

 \ Evidently, the matrix $ D = D_g $ is non-negative definite,  write $  D = D_g \ge \ge 0.  $ \par

\vspace{3mm}

 \ {\bf Definition 2.5.}\par

\vspace{3mm}

  \ Let $  V, \ V \subset R^d $ be open convex  central symmetric:  $ \forall x  \in V \ \Rightarrow -x \in V  $
 subset of whole space $  R^d $ containing some non-empty neighborhood of origin. We will denote the  collection
 of all such a sets  by $  S = S(R^d).  $\par

   \ We will denote by $  Y(V), \ V \in S(R^d) $ the set of all such a functions from the definition  2.4 which are defined only
on the set $  V. $ For instance, they are even, twice  continuous differentiable, convex, non-negative, is equal to zero
only at the origin  and $  \det D_g > 0,  $

$$
\lim_{x \to \partial V -0} \frac{\partial^d g }{\prod_{k=1}^d  \partial x_k} = \infty. \eqno(2.7)
$$

 \ Notation: $  V = \supp g.  $ \par

\vspace{4mm}

\section{ Definition and Banach structure of offered spaces. } \par

\vspace{4mm}

 \ {\bf  Definition 3.1. } \par

\vspace{3mm}

 \ Let the set $ V  $ be from the set
  $ S(R^d): V \in S(R^d)  $ and let the Young function  $  \phi(\cdot)  $ be from the set $  Y(V): \supp \phi = V. $\par

 \ We will say by definition likewise the one-dimensional case  that the {\it centered} random vector (r.v)
 $ \xi = \xi(\omega) = \vec{\xi} = (\xi(1), \xi(2), \ldots, \xi(d)) $ with values in the space $  R^d $ belongs to the
 space $ B_V(\phi), $  write $ \vec{\xi}\in  B_V(\phi),  $  if there exists certain non-negative constant $ \tau \ge 0 $ such that

$$
\forall \lambda \in V \ \Rightarrow
\max_{\vec{\epsilon}} {\bf E} \exp \left( \sum_{j=1}^d \epsilon(j) \lambda(j) \xi(j) \right) \le
\exp[ \phi(\lambda \cdot \tau) ]. \eqno(3.1)
$$

 \ The minimal value $ \tau $ satisfying (3.1) for all the values $  \lambda \in V, $
is named by definition as a $ B_V(\phi) \ $ norm of the vector $ \xi, $ write

$$
||\xi||B_V(\phi)  \stackrel{def}{=}
$$

 $$
 \inf \left\{ \tau, \ \tau > 0: \ \forall \lambda:  \ \lambda \in V \ \Rightarrow
 \max_{\vec{\epsilon}}{\bf E}\exp \left( \sum_{j=1}^d \epsilon(j) \lambda(j) \xi(j) \right) \le
 \exp(\phi(\lambda \cdot \tau)) \right\}. \eqno(3.2)
 $$

\vspace{3mm}

 \ For example, the {\it generating function } $  \phi_{\xi}(\lambda) $ for these spaces may be picked by the
 following natural way:

$$
\exp[ \phi_{\xi}(\lambda ) ] \stackrel{def}{=}
\max_{\vec{\epsilon} \in \Theta} {\bf E} \exp \left( \sum_{j=1}^d \epsilon(j) \lambda(j) \xi(j) \right),
 \eqno(3.2a)
$$
if of course the random vector  $  \xi $ is centered and has an exponential tail of distribution. This imply  that
the natural function  $ \phi_{\xi}(\lambda ) $ is finite on some non-trivial central symmetry neighborhood of origin,
or equivalently  the mean zero random vector  $  \xi $ satisfies the multivariate Kramer's  condition. \par

 \ Obviously, for the natural function $ \phi_{\xi}(\lambda )  $

$$
||\xi||B(\phi_{\xi}) = 1.
$$

 \ It is easily to see that this  choice of the generating function $ \phi_{\xi} $ is optimal, but in the practical using
often this function can not be calculated in explicit view, but there is a possibility to estimate its. \par

 \vspace{3mm}

 \ We agree to  take in the case when   $  V = R^d \hspace{4mm}  \ B_{R^d}(\phi):= B(\phi).  $ \par

 \ Note that the expression for the norm $ ||\xi||B_V(\phi)  $ dependent aside from the function $  \phi  $ and the set $ V, $
only on the distribution $  \Law(\xi). $ Thus, this norm and correspondent space  $  B(\phi)  $ are rearrangement invariant (symmetrical)
in the terminology of the classical book  \cite{Bennet1}, see chapters 1,2. \par

\vspace{4mm}

 \ {\bf Theorem 3.1.}  {\it  The space $ B_V(\phi) $ with respect to the norm $ || \cdot ||B_V(\phi) $ and
ordinary algebraic operations is a rearrangement invariant vector Banach space.} \par

 \vspace{3mm}

{\bf Proof.} Let us prove at first the {\it triangle inequality.}  We will denote for simplicity
$  ||\xi|| =   ||\xi||B_V(\phi)   $ and correspondingly  $  ||\eta|| =   ||\eta||B_V(\phi).   $ It is reasonable
to suppose $ \lambda \ge 0, \  0 < ||\xi||, \ ||\eta|| < \infty.  $ \par
 \ Assume for definiteness $   \vec{\epsilon} =\vec{1} = (1,1,\ldots,1). $ We conclude on the basis of  the direct definition of the
$  B_V(\phi) $ norm   for all the values $  \lambda \in V $

$$
{\bf E}  e^{ (\lambda, \ \xi)  } \le e^{\phi( ||\xi||  \cdot \lambda) }, \hspace{6mm}
{\bf E}  e^{ (\lambda, \eta)  } \le e^{\phi( ||\eta|| \cdot \lambda) }.
$$

 \ We apply the H\"older's inequality

$$
{\bf E} e^{(\lambda, \xi + \eta) } \le \sqrt[p]{ e^{ \phi(p \ ||\xi|| \ \lambda) }  } \cdot
\sqrt[q]{ e^{ \phi(q \ ||\eta|| \ \lambda)}},  \eqno(3.3)
$$
where as usually $ p,q = \const > 1, \ 1/p + 1/q = 1. $ We can choose in the estimate (3.3)

$$
p = \frac{||\xi|| + ||\eta||}{||\xi||}, \hspace{5mm} q = \frac{||\xi|| + ||\eta||}{||\eta||},
$$
and obtain after substituting into (3.3)

$$
{\bf E} e^{(\lambda, \xi + \eta) } \le e^{ \phi(\lambda \cdot (||\xi|| + ||\eta||) ) }, \eqno(3.4)
$$
therefore $  ||\xi + \eta|| \le ||\xi|| + ||\eta||. $ \par

 \ The equality  $  ||-\xi|| = ||\xi||  $ follows from the properties of parity of the function $  \phi(\cdot), $
the equality $ ||\alpha \ \xi|| = \alpha \ ||\xi||, \ \alpha = \const > 0 $ follows from the direct definition (3.2)
of norm in these spaces after simple calculations.\par

 \ Finally, the {\it completeness}  of our  $ B_V(\phi)  $ spaces follows immediately from the
one-dimensional case as long as take place here  the  coordinatewise completeness. \par

\vspace{4mm}

 \ {\bf Remark 3.1.} In the article of  Buldygin V.V. and Kozachenko Yu. V. \cite{Buldygin1} was
 considered a particular case when $  V = R^d  $ and $  \phi(\lambda) = \phi^{(B)}(\lambda)  = 0.5(B\lambda,\lambda),  $ where
$  B  $  is non-degenerate positive definite symmetrical matrix, as a direct generalization of the one-dimensional one notion,
belonging to J.P.Kahane \cite{Kahane1}.  \par
 \ The correspondent random vector $ \vec{\xi}  $ was named in \cite{Buldygin1} as a subgaussian r.v. relative the
matrix $  B: $ \par

$$
\forall \lambda \in R^d \ \Rightarrow {\bf E} e^{(\lambda, \xi)} \le e^{0.5 (B \lambda, \lambda) \ ||\xi||^2 }.
$$

 \ We will write in this case  $  \xi \in \Sub(B)  $ or more precisely $ \Law \xi \in \Sub(B). $ \par

\vspace{4mm}

 \ {\bf Remark 3.2.} Let the function $  \phi(\cdot) \in Y(V) $ be a given. Define the following operation between two
non - negative numbers $  a $ and $  b: $

$$
a \odot b = a \odot_{\phi} b \stackrel{def}{=} \inf \{c, \ c > 0, \ \forall \lambda \in V \ \Rightarrow
\phi(c \cdot \lambda) \ge \phi(a \cdot \lambda)  + \phi(b \cdot \lambda) \}.
$$

 \ It is clear that $ 0 \odot a  = a, \ a > 0;  $

$$
 a \odot b = b \odot a,  \ (a \odot b) \odot c  = a \odot ( b \odot c);
$$

$$
  (\beta \cdot a) \odot (\beta \cdot b) = \beta \cdot  (a \odot b), \ \beta = \const \ge 0.
$$

 \vspace{4mm}

 \ {\it We propose:} if both the random vectors $ \xi, \eta  $ are independent and  belongs to the space $ B_V(\phi),  $ then

$$
||\xi + \eta|| \le ||\xi|| \odot_{\phi} ||\eta||.
$$

 \ Indeed, by virtue of independence

$$
{\bf E}e^{(\lambda, (\xi + \eta)) } = {\bf E}e^{(\lambda,\xi)  } \times  {\bf E}e^{(\lambda,\eta) } \le
$$

$$
e^{ \phi(||\xi|| \ \lambda)} \cdot e^{ \phi(||\eta|| \ \lambda)} \le e^{ \phi( (||\xi|| \odot ||\eta||) \ \lambda ) }.
$$

 \ As a slight consequence:  if $  \{\xi_i \}, \ i = 1,2,\ldots,n $  be a sequence of independent (centered) multidimensional
r.v.  to at the same space $  B(\phi).  $ Denote

$$
S(n) = n^{-1/2} \sum_{i=1}^n  \xi_i, \eqno (3.5)
$$
then

$$
 ||S(n)|| B(\phi) \le n^{-1/2}
 \left\{ \  ||\xi_1|| B(\phi) \odot_{\phi} ||\xi_2|| B(\phi) \odot_{\phi} \ldots  \odot_{\phi} ||\xi_n|| B(\phi) \ \right\}. \eqno(3.6)
$$

\vspace{4mm}

 \ {\bf Remark 3.3.}  Suppose the r.v. $  \xi  $ belongs to the space $  B(\phi);  $  then evidently
 $  {\bf E} \vec{\xi} = 0. $ \par
  \ Suppose in addition that it has there the unit norm, then

 $$
 \Var ( \vec{\xi}) \le \le   D_\phi; \hspace{4mm} \Leftrightarrow D_\phi \ge \ge \Var ( \vec{\xi}).
  \eqno(3.7)
 $$

\vspace{3mm}

 \ It is interest to note that there exist many random vectors $  \eta = \vec{\eta} $ for which

$$
{\bf E} e^{(\lambda,\eta)} \le e^{0.5 \ (\Var(\eta) \lambda, \ \lambda) }, \ \lambda \in R^d, \eqno(3.8)
$$
see e.g.  \cite{Buldygin2}, \cite{Buldygin3}, chapters 1,2; \cite{Ostrovsky1}, p.53. V.V.Buldygin and
Yu.V.Kozatchenko in \cite{Buldygin2} named these vectors {\it  strictly subgaussian;}  notation
$ \xi \in \SSub $ or equally $ \Law (\xi) \in \SSub. $ \par

 \ V.V.Buldygin and Yu.V.Kozatchenko found also some interest applications of these notions. \par

\vspace{4mm}

\section{Characterization features.}

\vspace{4mm}

 \ Statement of problem: given a function $  \phi(\cdot) $  from the set $ Y = Y(R^d). $ Question: under what
additional conditions there exists a mean zero random   vector $  \xi  $ which may be defined on the appropriate
probability space such that

$$
\forall \lambda \in R^d \ \Rightarrow {\bf E} e^{ (\lambda,\xi) } = e^{\phi(\lambda)}, \eqno(4.1)
$$
 or equally $ \phi(\lambda) =  \phi_{\xi}(\lambda)  $ for some random vector $  \xi. $ \par

  \ It is clear that we hinted at multivariate Bernstein's theorem, see \cite{Feller1}, \cite{Jozsa1}. Evidently, $  \phi(0) = 0.  $
  \par

\vspace{3mm}

 \ Denote by $  \mu_{\xi}(\cdot) $ the (Borelian) distribution of r.v. $  \xi:  \ \mu_{\xi}(A) := {\bf P} (\xi \in A), \ A \subset R^d, $
so that

$$
{\bf E} e^{ (\lambda,\xi) } = \int_{R^d} e^{ (\lambda, x) }  \ \mu_{\xi}(dx). \eqno(4.2)
$$

 \ It will be presumed that the measure $ \mu_{\xi} $ and some ones have exponential decreasing tails, so that the integral
 in (4.2) converges for all the values $ \lambda \in R^d. $ \par

 \ Let us consider the  following integral

$$
I_+(\lambda) := \int_{R^d_+} e^{ (\lambda, x) }  \ \mu_{\xi}(dx), \eqno(4.3)
$$
here and in the sequel

$$
R^d_+ \stackrel{def}{=}  \{ \vec{x}: \ \forall j \ \Rightarrow x(j) \ge 0 \}.
$$

 \ We deduce

$$
\frac{\partial^k I_+(\lambda)}{ \partial \lambda(1)^{k(1)} \ \partial \lambda(2)^{k(2)} \  \ldots \partial \lambda(d)^{k(d) }} =
\int_{R^d_+}  e^{ (\lambda, x) } \cdot \prod_{j=1}^d x(j)^{k(j)} \cdot \mu_{\xi}(dx),
$$
therefore

$$
\frac{\partial^{ \vec{k}} \ I_+(\lambda) }{ \partial \vec{\lambda}^{{\vec k}}} =
\frac{\partial^{|k|} I_+(\lambda)}{ \partial \lambda(1)^{k(1)} \ \partial \lambda(2)^{k(2)} \  \ldots \partial \lambda(d)^{k(d) }}
 \ge 0, \eqno(4.4)
$$
where of course $  k, k(j) = 0,1,2, \ldots $ and $ \sum_j k(j) = k. $ \par

 \ The function satisfying the inequalities (4.4) for all the integer non - negative vectors $ \vec{k}  $
 are named  {\it  absolutely monotonic. } \par

 \ The multivariate  version of the classical Bernstein's theorem, see for example  \cite{Jozsa1}, tell us that if the function,
 say $ I(\lambda)   $ is infinitely many times differentiable and satisfies the inequality (4.4) for all the integer positive vectors
 $ \vec{k}, $ then there exists a Borelian {\it measure} $ \nu $ on the set $ R^d_+  $ such that

$$
I(\lambda) = \int_{R^d_+} e^{ (\lambda, x) }  \ \nu(dx). \eqno(4.5)
$$

\vspace{3mm}

 \ Let $ \vec{\epsilon} = \epsilon = \{ \epsilon(1),   \epsilon(2), \ldots,   \epsilon(d) \} $ be non-random  element of the set $  \Theta. $
Recall that $ \vec{\epsilon}  $ is $  d \ - $ dimensional  numerical vector, $  d = 2,3,\ldots, $ whose components take the values
$  \pm 1 $ only. We associate for all such a vector $ \vec{\epsilon} $ the following octant in the whole space $  R^d  $

$$
Z(\vec{\epsilon}) = Z(\epsilon) = \{ \vec{x}: \ \forall j \ \Rightarrow  \epsilon(j) \ x(j) \ge 0 \} \eqno(4.6)
$$
and correspondent integral

$$
I(Z(\epsilon),\lambda) = \int_{Z(\epsilon)} e^{ (\lambda, x) }  \ \nu(dx). \eqno(4.6a)
$$

 \ For instance,
$$
Z(\vec{1}) = R^d_+; \hspace{5mm} I(Z(\vec{1}),\lambda) = \int_{R^d_+} e^{ (\lambda, x) }  \ \nu(dx) = I(\lambda).
$$

 \ We get alike (4.4)

$$
\prod_{j=1}^d \epsilon(j)^{k(j)} =
\sign \frac{\partial^{ \vec{k}} \ I(Z(\vec{\epsilon}),\lambda) }{ \partial \vec{\lambda}^{{\vec k}}} =
$$

$$
\sign \frac{\partial^{|k|} I(Z(\vec{\epsilon}),\lambda)}{ \partial \lambda(1)^{k(1)} \
\partial \lambda(2)^{k(2)} \  \ldots \partial \lambda(d)^{k(d) }}=
\vec{\epsilon}^{\vec{k} }.  \eqno(4.7)
$$

\vspace{3mm}

 \ {\bf Definition 4.1.} Let the vector $ \vec{\epsilon} \in \Theta  $ and related octant $ Z(\vec{\epsilon})  $
 be a given. The numerical function  $  F = F(\lambda), \ \lambda \in R^d $ belongs  by definition to the class
 $  K( \vec{\epsilon} ) = K(\epsilon),  $  or on the other words, is monotonic relative the octant $  Z(\vec{\epsilon}),  $
 iff it satisfied the restriction (4.7) for all the  values $ \vec{k} : $

$$
\forall \lambda \in R^d \ \Rightarrow
\sign \frac{\partial^{ \vec{k}} \ F(\lambda) }{ \partial \vec{\lambda}^{{\vec k}}}=
\vec{\epsilon}^{\vec{k} }.  \eqno(4.8)
$$

 \ It follows immediately from the multivariate Bernstein's theorem after replacing  $  x \to \epsilon \otimes x $
the following fact. \par

\vspace{3mm}

 \ {\bf Lemma 4.1.} The function $ F = F(\lambda), \ \lambda \in R^d $ belongs   to the class   $ K( \vec{\epsilon} ) = K(\epsilon)  $
 iff it has a representation

$$
F(\lambda) = \int_{Z(\epsilon)} e^{(\lambda, x)} \ \zeta(dx) \eqno(4.9)
$$
for some finite Borelian measure $ \zeta(\cdot). $  \par

\vspace{4mm}

 {\bf Theorem 4.1.} The infinite differentiable
  function $ \phi(\lambda)  $  such that $ \phi(0) = 0  $
  may be represented  on the form (4.1): $  \phi = \phi_{\xi}  $
for some random vector $ \ \xi = \vec{\xi} \ $ defined on some sufficiently rich probability space,
iff the function $  \exp \phi(\lambda) $ has a representation of the form

$$
 \exp \phi(\lambda)= \sum_{\vec{\epsilon}  \in \Theta} F_{ \vec{\epsilon}}(\lambda), \eqno(4.10)
$$
where each function  $ F_{ \vec{\epsilon}}(\lambda) $ is monotonic relative the octant $  Z(\vec{\epsilon}): \
F_{ \vec{\epsilon}}(\cdot) \in K( \vec{\epsilon}) $  and

$$
\sum_{\vec{\epsilon}  \in \Theta} F_{ \vec{\epsilon}}(0) = 1.\eqno(4.10a)
$$

\vspace{4mm}

 \ {\bf Proof.} Assume first of all that the expression (4.10) take place. Denote by $  \chi_A(x)  $ the indicator function
of the measurable set $ A. $  We deduce using the proposition of lemma 4.1:

$$
\sum_{\vec{\epsilon}  \in \Theta} F_{ \vec{\epsilon}}(\lambda) =
\sum_{\vec{\epsilon}  \in \Theta} \int_{Z(\epsilon)} \exp(\lambda,x) \  \nu_{\epsilon}(dx) =
\sum_{\vec{\epsilon}  \in \Theta} \int_{R^d} \chi_{Z(\epsilon)}(x) \exp(\lambda,x)  \ \nu_{\epsilon}(dx) =
$$

$$
 \int_{R^d} \exp(\lambda,x) \sum_{\vec{\epsilon}  \in \Theta} \ \chi_{Z(\epsilon)}(x) \ \nu_{\epsilon}(dx) =
 \int_{R^d} \exp(\lambda,x) \ \nu(dx), \eqno(4.11)
$$
where

$$
\nu(A) = \int_A \sum_{\vec{\epsilon}  \in \Theta} \ \chi_{Z(\epsilon)}(x) \ \nu_{\epsilon}(dx). \eqno(4.11a)
$$

 \  Substituting $  \lambda = 0 $ into (4.11), we conclude that $  \nu $ is actually a probabilistic measure and

$$
 \exp \phi(\lambda)= \int_{R^d} \exp(\lambda,x) \nu(dx) = {\bf E} \exp(\lambda, \gamma),
 \eqno(4.12)
$$
where the r.v. $  \gamma $ has the distribution $  \nu.  $\par
 \ Conversely, let

$$
e^{\phi(\lambda)} =  {\bf E} e^{(\lambda, \xi)}
$$
for some $  d \ - $ dimensional random vector $ \xi,  $ then

$$
e^{\phi(\lambda)} =  \int_{R^d} e^{(\lambda, \xi)} \ \mu_{\xi}(dx) =
$$

$$
\sum_{\vec{\epsilon}  \in \Theta} \int_{Z(\epsilon)} \exp(\lambda,x) \  \mu_{\epsilon}(dx) =
\sum_{\vec{\epsilon} \in \Theta}  F_{ \vec{\epsilon}  }(\lambda),
$$
where the functions  $ F_{ \vec{\epsilon}  }(\lambda): $

$$
F_{ \vec{\epsilon}}(\lambda)=   \int_{Z(\epsilon)} \exp(\lambda,x) \  \mu_{\epsilon}(dx)
$$
satisfy the conditions of lemma 4.1.\par

\vspace{4mm}

\section{ Linear transforms.}

\vspace{4mm}

 \ Statement of problem: given a $  d \ - $ dimensional (centered)  random vector $ \xi $ from the space   $  B(\phi), $ in
particular

$$
{\bf E} e^{(\lambda,\xi)} \le e^{ \phi(||\xi|| \cdot \lambda) }, \ \phi \in Y(R^d), \eqno(5.1)
$$
and let $  A: R^d \to R^d $ be linear operator acting from the space $  R^d $  into itself.
It is required  to estimate  some $  B(\tilde{\phi}) $ norm of the linear transformed random vector $ \eta = A\cdot \xi.  $\par

 \ Evidently, the r.v. $ \eta $ belongs to at the same space $  B(\phi),  $ as long as it may be represented as linear combination
of the elements of these spaces. We intend to obtain more exact estimate. \par

 \ Define the following function from at the same space $  \Phi: $

$$
\phi_A(\vec{\lambda}) := \phi(A^* \  \vec{\lambda} ).\eqno(5.2)
$$

\vspace{4mm}

 \ {\bf Theorem 5.1.}

$$
|| A \xi||B(\phi) = ||\xi||B(\phi_A). \eqno(5.3)
$$

\vspace{4mm}

 {\bf Proof } is very simple: $ {\bf E} \exp(\lambda, A \xi) =  $

$$
 {\bf E} \exp( A^*\lambda, \xi) \le  \exp( \phi(||\xi|| \cdot A^* \lambda )) = \exp (\phi_A(||\xi||)),\eqno(5.4)
$$
therefore

$$
|| A \xi||B(\phi) \le ||\xi||B(\phi_A).
$$
 Inverse  inequality may be justify analogously. \par

 \vspace{4mm}

 \ {\bf Example 5.1.} \\

 \ Let $  \xi $  be subgaussian random vector relative  certain positive definite symmetrical  matrix $  R. $
Then

$$
|| A \xi|| \Sub(R) = ||\xi||\Sub (A \ R \ A^*). \eqno(5.5)
$$

\vspace{4mm}

 \ We introduce now  some generalization of the well - known notion of a $  \Delta_2 $ condition on the
multidimensional case. \par

 \vspace{4mm}

 \ {\bf Definition 5.1.} We will say that the function $  \phi = \phi(\lambda), \lambda \in R^d $ satisfies the
multidimensional $ \vec{\Delta}_2  $  condition, write $ \phi(\cdot) \in \vec{\Delta}_2, $
 if for arbitrary matrix $  A: R^d  \to R^d $ there exists a
non - negative constant, which we will denote by $  |||A|||_{\phi}, $ such that

$$
\forall \lambda \in R^d \ \Rightarrow \phi(A^* \lambda) \le \phi( |||A|||^2_{\phi} \cdot \lambda). \eqno(5.6)
$$

 \ It is easily to verify that the functional $  \ A \to |||A|||_{\phi}  $ is actually certain matrix norm. \par

 \ It follows from proposition of the example (5.1) that the quadratic form $ \lambda \to (R\lambda,\lambda)  $ satisfies the
 $ \vec{\Delta}_2  $  condition.\par

\vspace{4mm}

 \ {\bf Proposition 5.1.} \ It follows  immediately from the relation (5.4) that if the function $ \phi \in Y(R^d) $
belongs also the class $ \vec{\Delta}_2, $ then under at the same notations and restrictions

$$
|| A \xi||B(\phi) \le  |||A|||_{\phi} \cdot  ||\xi||B(\phi). \eqno(5.7)
$$

\vspace{4mm}

\section{ Tails behavior. } \par

\vspace{4mm}

{\bf A. Upper estimate.} \par

\vspace{3mm}

 \ Let $  \phi = \phi(\lambda), \ \lambda \in V \subset R^d $ be arbitrary non-negative real valued function, as in the definition 2.3,
which is finite on some-empty neighborhood of origin.  Suppose for given centered $  d  \ - $ dimensional random vector $ \xi = \vec{\xi}  $

$$
{\bf E} e^{(\lambda,\xi)} \le e^{ \phi( \lambda) }, \ \lambda \in V.  \eqno(6.1)
$$

 \ On the other words, $ ||\xi||B(\phi) \le 1. $ \par

\vspace{4mm}

 \ {\bf Proposition 6.1.}  For all non-negative vector $ x = \vec{x}  $ there holds

 $$
 U(\vec{\xi}, \vec{x}) \le \exp \left( - \phi^*(\vec{x}) \right) \ - \eqno(6.2)
 $$
the multidimensional generalization of Chernov's inequality. \par

\vspace{3mm}

 \ {\bf Proof.}  Let for definiteness $  \ x_j > 0; $ the case when $ \exists k \ \Rightarrow x_k \le 0  $ may be
considered analogously. \par
 \ We use the ordinary Tchebychev's  inequality

$$
U(\vec{\xi}, \vec{x}) \le  \frac{ e^{ \phi( \lambda) }}{ e^{( \lambda,x )}  } =
e^{ -( \lambda,x ) + \phi( \lambda)}, \ \lambda > 0. \eqno(6.3)
$$

 \ Since the last inequality  (6.3) is true for arbitrary non - negative vector $  \lambda \in V,  $

$$
U(\vec{\xi}, \vec{x}) \le \inf_{\lambda \in V} e^{( \lambda,x ) - \phi( \lambda)} =
$$

$$
\exp \left\{- \sup_{\lambda \in V} [ \ ( \lambda,x ) - \phi( \lambda) \ ] \right\}  = \exp \left( - \phi^*(\vec{x}) \right). \eqno(6.4)
$$

\vspace{3mm}

{\bf B. Lower estimate.} \par

\vspace{3mm}

 \ {\bf Theorem 6.1.} Let the function $  \phi(\cdot)  $ be from the set $  S(R^d),  $  (Definition 2.5.)
 Suppose the mean zero random vector $  \xi = \vec{\xi} $ satisfies the condition (6.2) for all non-negative deterministic
 vector $  \vec{x}: \ \forall j = 1,2,\ldots,d \ x_j > 0  $

 $$
 U(\vec{\xi}, \vec{x}) \le \exp \left( - \phi^*(\vec{x}) \right).    \eqno(6.5)
 $$
 \ We propose that r.v. $ \ \vec{\xi} \ $ belongs to the space $  B(\phi): \ \exists  \ C = C(\phi) \in (0,\infty), $

$$
{\bf E} e^{(\lambda,\xi)} \le e^{ \phi( C \cdot \lambda) }, \ \lambda \in R^d. \eqno(6.6)
$$

 \vspace{3mm}

 \ {\bf Proof}  is likewise to the one-dimensional case, see  \cite{Kozachenko1}, \cite{Ostrovsky1}, p. 19-40;
 see also \cite{Ostrovsky5}.\par
 \ Note first of all that the estimate (6.6) is obviously satisfied   for the values $  \lambda = \vec{\lambda} $
 from the Euclidean unit ball of the space $   R^d: \ |\lambda| \le 1, $ since the r.v. $  \xi $ is centered and has a very
 light (exponential decreasing) tail of distribution. It remains to consider further only the case when $  |\lambda| \ge 1.$ \par

\vspace{3mm}

 \ We have using integration  by parts

$$
 {\bf E} e^{(\lambda,\xi)}  \le \prod_{k=1}^d |\lambda_k| \cdot \int_{R^d} e^{(\lambda,x) - \phi^*(x) } \ dx \stackrel{def}{=}
  \prod_{k=1}^d |\lambda_k| \cdot  I_{R^d}(\lambda).
$$

 \ It is sufficient to investigate the main part of the last integral, indeed

$$
I_{R_+^d}(\lambda) := \int_{R_+^d} e^{(\lambda,x) - \phi^*(x) } \ dx. \eqno(6.7)
$$

 \ Let $  \gamma = \const \in (0,1).  $   We have using Young inequality

$$
(\lambda,x) = (\gamma \lambda, x/\gamma) \le \phi^*(\gamma x) +   \phi^{**}(\lambda/\gamma),
$$
 and after substituting into (6.7)

$$
I_{R_+^d}(\lambda) \le  e^{ \phi^{**}(\lambda/\gamma)} \cdot \int_{R_+^d} e^{\phi^*(\gamma x) - \phi^*(x) } \ dx =
C_2(\gamma,\phi) \cdot e^{ \phi^{**}(\lambda/\gamma)}. \eqno(6.8)
$$

  \ We conclude by virtue of theorem Fenchel -Moraux $ \phi^{**}(\lambda) = \phi(\lambda),  $ therefore

$$
 {\bf E} e^{(\lambda,\xi)}  \le C_3(\gamma,\phi) \ \left[ C_4 (\gamma,\phi) + C_5 (\gamma,\phi) \prod_{k=1}^d |\lambda_k| \ \right] \cdot
e^{ \phi(\lambda/\gamma)} \le
$$

$$
 e^{ \phi(\lambda/\gamma_2)}, \ |\lambda| \ge 1, \eqno(6.9)
$$
 where $  \gamma_2 = \gamma_2(\phi) = \const \in (0,1). $ \par

 \  Another details are simple and may be omitted. \par

\vspace{4mm}

 \ As a slight consequence:\par

\vspace{4mm}

 \ {\bf  Corollary 6.1. }   Let as above the function $  \phi(\cdot)  $ be from the set $  S(R^d),  $  see the definition 2.5.
The centered non-zero random vector $  \xi  $ belongs to the space  $   B(\phi): $

$$
\exists C_1 \in (0,\infty), \ \forall \lambda \in R^d
 \Rightarrow  {\bf E} e^{(\lambda,\xi)} \le e^{ \phi( C_1 \cdot \lambda) }, \ \lambda \in R^d
$$
if and only if

$$
\exists C_2 \in (0,\infty), \ \forall \ x \in R^d_+  \Rightarrow \  U(\vec{\xi}, \vec{x}) \le \exp \left( - \phi^*(\vec{x}/C_2) \right).
$$

 \ More precisely, the following implication holds: there is finite positive constant $  C_3 =  C_3(\phi)  $ such that
  for arbitrary non-zero centered r.v. $  \xi: \ ||\xi|| = ||\xi||B(\phi)  < \infty \ \Leftrightarrow $

$$
 \forall \lambda \in R^d
 \Rightarrow  {\bf E} e^{(\lambda,\xi)} \le e^{ \phi( ||\xi|| \cdot \lambda) }
$$
iff

$$
\exists C_3(\phi) \in (0,\infty) \ \forall \ x \in R^d_+  \Rightarrow \  U(\vec{\xi}, \vec{x}) \le
\exp \left( - \phi^*(\vec{x}/(C_3 / ||\xi||) \right).
\eqno(6.10)
$$

\vspace{3mm}

\ {\bf  Corollary 6.2. } Assume the non-zero centered random vector $  \xi = (\xi(1), \xi(2), \ldots, \xi(d)) $ belongs
to the space $  B(\phi):  $

$$
{\bf E} e^{(\lambda,\xi)} \le e^{ \phi(||\xi|| \cdot \lambda) }, \ \phi \in Y(R^d),
\eqno(6.11)
$$
and let $  y  $ be arbitrary positive non-random number. Then $  \forall y  > 0 \ \Rightarrow  $

$$
{\bf P} \left( \min_{j = 1,2,\ldots,n} |\xi(j)| > y  \right) \le 2^d \cdot
\exp \left(- \phi^*(y/||\xi||,y/||\xi||, \ldots, y/||\xi||) \right). \eqno(6.12)
$$

\vspace{4mm}

{\bf Example 6.1.} Let as before $  V = R^d  $ and $  \phi(\lambda) = \phi^{(B)}(\lambda)  = 0.5(B\lambda,\lambda),  $ where
$  B  $  is non-degenerate positive definite symmetrical matrix, in particular $  \det B > 0. $ It follows from theorem 6.1
that the (centered) random vector $ \xi $ is subgaussian relative the matrix $ B: $

$$
\forall \lambda \in R^d \ \Rightarrow {\bf E} e^{ (\lambda, \xi)} \le e^{0.5 (B \lambda, \lambda) ||\xi||^2 }.
$$
iff for some finite positive constant $  K = K(B,d) $ and for any  non-random positive vector $  x = \vec{x} $

$$
U(\xi,x) \le e^{- 0.5 \ \left( (B^{-1}x,x)/(K||\xi||^2) \right) }. \eqno(6.13)
$$

\vspace{4mm}

\section{ Relation with moments. } \par

\vspace{4mm}

 \ We intend in this section to simplify the known proof the moment estimates (1.6)-(1.6a) for the one-dimensional r.v. and
extend  obtained result on the multivariate case. \par

 \ We will use the following elementary inequality

$$
x^r \le \left( \frac{r}{\lambda \ e} \right)^r \cdot e^{\lambda \ x}, \ r,\lambda,x > 0, \eqno(7.0)
$$
and hence

$$
|x|^r \le \left( \frac{r}{\lambda \ e} \right)^r \ \cosh (\lambda x), \  \ r,\lambda > 0, \ x \in R. \eqno(7.0a)
$$

 \ As a consequence: let $  \xi $ be non-zero  one-dimensional mean zero random variable belonging to the space $  B(\phi).  $
Then

$$
{\bf E} |\xi|^r \le 2 \  \left( \frac{r}{\lambda \ e} \right)^r  \ e^{ \phi(\lambda ||\xi||) }, \ \lambda > 0. \eqno(7.1)
$$

  \ Author of the works \cite{Kozachenko1}, \cite{Ostrovsky1}, chapters 1,2 choose  in the inequality (7.1) the value

$$
\lambda = \lambda_0 := \phi^{-1}(r/||\xi||)
$$
and obtained the relations (1.6), (1.6a). We intend here to choose the value $ \lambda $ for reasons of optimality. \par
 \ We hope that this method is more simple and allows easy multivariate generalization. \par

 \ In detail, introduce the function

$$
\Phi(\mu) := \phi \left(e^{\mu} \right), \ \mu \in R. \eqno(7.2)
$$

 \ Let temporarily  for simplicity in (7.1) be $ ||\xi|| = ||\xi||B(\phi) = 1. $ One can rewrite  (7.1) as follows.

$$
{\bf E} |\xi|^r \le 2 \ r^r \ e^{-r} \ e^{ - r \ln \lambda + \phi(\lambda ) } =
$$

$$
 2 \ r^r \ e^{-r} \ e^{ - r \mu + \phi \left(e^{\mu} \right) } = 2 \ r^r \ e^{-r} \ e^{ -\left(r \mu - \Phi(\mu) \right) }, \eqno(7.3)
$$
and we deduce after minimization over $ \mu \ $ or $ (\lambda)  $

$$
{\bf E} |\xi|^r \le  2 \ r^r \ e^{-r} \ e^{- \Phi^*(r) }.
$$

 \ So, we proved in fact the following statement. \par

 \vspace{3mm}

{\bf Proposition 7.1.} Let $ \phi (\cdot) $  be arbitrary non - negative continuous function and let the centered numerical
r.v. $  \xi $  be such that $  \xi \in B(\phi) $ or equally

$$
U(\xi,x) \le \exp(-\phi^*(x)), x \ge 0.
$$

Then

$$
|\xi|_r \le 2^{1/r} \ r \ e^{-1} \  e^{- \Phi^*(r)/r } \ ||\xi||B(\phi), \ r > 0. \eqno(7.4)
$$

\vspace{3mm}

 \  We want investigate now the inverse conclusion. Namely, let the mean zero numerical  r.v. $ \xi  $ be such that

$$
|\xi|_r \le C_0 \ r \ e^{- \Phi^*(r)/r }.  \eqno(7.5)
$$

\vspace{3mm}

 Let the function $  \Phi(\cdot)  $ be from the set $  S(R) $ and suppose the real valued centered
r.v.  $  \xi $ satisfies the estimate (7.5) at least for the integer even  values $  r = 2 m, \ m = 1,2, \ldots: $

$$
|\xi|_{2m} \le  K \cdot (2m) \cdot e^{- \Phi^*(2m)/(2m)}, \ K = \const \in (0,\infty).  \eqno(7.6)
$$

 \ The last inequality may be rewritten on the language of the $  G\psi  $ spaces as follows. Define the correspondent
 $  \psi \ -  $ function

$$
\psi(m) = \psi_{\phi}(m) = (2m) \cdot e^{- \Phi^*(2m)/(2m)}, \ m = 1,2, \ldots;
$$
then the constant $  K  $ in (7.6) has a form

$$
K = ||\xi||G\psi_{\phi} = \sup_m \left[ \frac{|\xi|_{2m}}{\psi_{\phi}(m)}\right]. \eqno(7.7)
$$

\vspace{3mm}

{\bf Theorem 7.1.} Suppose $  \phi \in S(R^1). $  Let $ {\bf E} \xi =0  $ and
$ K = ||\xi||G\psi_{\phi} < \infty. $ Then $ \xi \in B(\phi)  $ and of course

$$
|| \xi||B(\phi) \le C_1 \ ||\xi||G\psi_{\phi} \le  C_2 ||\xi||B(\phi), \ 0 < C_1 = \const \le C_2  < \infty.
\eqno(7.8)
$$

 \vspace{3mm}

 \ {\bf Proof.} It remains to prove only the left-hand side of the bilateral inequality (7.8).
 We can assume without loss of generality $  {\bf E} \xi = 0 $ and

$$
|\xi|_{2m} \le  (2m) \cdot e^{- \Phi^*(2m)/(2m)}.  \eqno(7.9)
$$

 \ Let us consider the so - called moment generating function $ g(\lambda) = g_{\xi}(\lambda) $ for the r.v. $  \xi: $

$$
g(\lambda) = g_{\xi}(\lambda) \stackrel{def}{=} {\bf E} e^{\lambda \xi}. \eqno(7.10)
$$

 \ Our target here is equivalent to justify the key estimate: $  \exists = C(m) \in (0,\infty)  $ such that

$$
 g_{\xi}(\lambda) = {\bf E} e^{\lambda \xi} \le e^{\phi(C(m) \lambda)}, \ \lambda \in R. \eqno(7.11)
$$

 \ The last inequality (7.11) is evident for the "small" values $  \lambda, $ say for $  |\lambda | \le 1. $\par

 \vspace{3mm}

 \ Let now $  |\lambda | \ge 1,  $ for definiteness let $ \lambda \ge 1.  $ We deduce taking into account the  centeredness
of the r.v. $  \xi: $

$$
{\bf E} \cosh{\lambda \xi} = 1 + \sum_{m=1}^{\infty} \frac{\lambda^{2m}}{ (2m)!}  \cdot {\bf E} \xi^{2m} \le
$$

$$
1 + \sum_{m=1}^{\infty} \frac{ C^m \ \lambda^{2m}}{ (2m)!} \cdot (2m)^{2m} \cdot e^{ - \Phi^*(2m) }  \le
1 + \sum_{m=1}^{\infty} \ C_2^m \cdot \lambda^{2m} \cdot e^{ - \Phi^*(2m) }.
$$
 \ We apply the so-called discrete analog of the saddle - point method, see
\cite{Sachkov1}, chapter 3, sections 3, 4:

$$
{\bf E} \cosh{\lambda \xi} \le \exp \left( \sup_{z \ge 1} \left[ z \cdot \ln ( C_3\lambda) - \Phi^*(z) \right]  \right) =
$$

$$
\exp \left(  \Phi^{**} ( \ln C_4 \lambda)   \right) = \exp \left(  \Phi ( \ln C_4 \lambda)   \right)= \exp \phi(C_4\lambda)
\eqno(7.12)
$$
 again by virtue of theorem  Fenchel-Moraux. \par

\vspace{3mm}

 \ This completes the proof of theorem 7.1. \par

\vspace{3mm}

 {\bf Example 7.1.}  Suppose that the function $ \phi (\cdot)  $  is from the set $  S(R^1) $ be such that for some
 constant $  p  > 1 $

$$
\phi(\lambda) = \phi_p(\lambda) \le C_1 \ |\lambda|^p, \ |\lambda| >1.
$$
 \ Let  also the centered non-zero random variable  $  \xi $ belongs to the space $ G\psi_p.  $
Then

$$
|\xi|_r \le C_2(p) \ r^{1/q} \ ||\xi||B(\phi_r), \ r \ge 1, \ q = p/(p-1),
$$
and the inverse conclusion is also true: if $ {\bf E} \xi = 0 $ and if for some constant $  K $

$$
\forall r \ge 1 \ \Rightarrow  |\xi|_r \le K \ r^{1/q},  \ q = p/(p-1),
$$
then $ \xi(\cdot) \in B(\phi_p) $  and wherein $ ||\xi||  B(\phi_p) \le C_3 \ K.   $ \par

\vspace{4mm}

 {\bf Example 7.2.}  Suppose now that the function $ \phi (\cdot)  $ is from the set $  S(-K, K), \ 0 < K= \const < \infty $
 be such that

$$
\phi(\lambda) = \phi^{(K)}(\lambda)  \le \frac{ C_4}{K - \ |\lambda|}, \ |\lambda| < K.
$$
 \ Let  also the centered non - zero random variable  $  \xi $ belongs to this space $ G\psi^{(K)}.  $
Then

$$
|\xi|_r \le C_5(K) \ r \ ||\xi||B(\phi^{(K)}),
$$
and likewise the inverse conclusion is also true: if $ {\bf E} \xi = 0 $ and if for some finite positive constant $  K $

$$
\forall r \ge 1 \ \Rightarrow  |\xi|_r \le K \ r,
$$
then $ \xi(\cdot) \in B(\phi^{(K)}) $  and herewith $ ||\xi|| B(\phi^{(K)}) \le C_6 \ K.   $ \par

\vspace{4mm}

 \ We need  getting to the presentation of the multidimensional case to extend our notations and restrictions.
 In what follows in this section the variables $ \lambda, r, x, \xi  $  are as before vectors from the space $  R^d, \ d = 2,3,\ldots, $
and besides $  r = \vec{r} = \{ r(1), r(2), \ldots, r(d) \}, \ r(j) \ge 1.  $\par

 \ Vector notations:

 $$
  |r| = |\vec{r}| = \sum_j r(j), \hspace{4mm} |\xi| = | \ \vec{\xi} \ | = \{ |\xi(1)|, |\xi(2)|, \ldots, |\xi(d)| \} \in R^d_+,
 $$

$$
 \vec{x} \ge 0 \ \Leftrightarrow \forall j \hspace{3mm} x(j) \ge 0;
$$

 $$
 x^r = \vec{x}^{\vec{r}} = \prod_{j=1}^d x(j)^{r(j)}, \ \vec{x} \ge 0,
 $$

$$
\ln \vec{\lambda} = \{  \ln \lambda(1), \ \ln \lambda(2), \ldots, \ \ln \lambda(d) \}, \hspace{3mm} \vec{\lambda} > 0,
$$

$$
e^{\vec{\mu}} = \{ e^{\mu(1)}, \ e^{\mu(2)}, \ldots, \  e^{\mu(d)} \},
$$

$$
\Phi(\mu) = \Phi(\vec{\mu}) = \phi \left(e^{\vec{\mu}} \right),
$$

 $$
 \frac{r}{\lambda \cdot  e}  =  \frac{\vec{r}}{ \vec{\lambda} \ e} = \prod_{j=1}^d \left( \frac{r(j)}{ e \lambda(j)} \right) =
 e^{-|r|} \cdot \prod_{j=1}^d \left( \frac{r(j)}{\lambda(j)} \right),
 $$

$$
|\xi|_r = |\vec{\xi}|_{\vec{r}} = \left( {\bf E}|\vec{\xi}|^{\vec{r}} \right)^{1/|r|}.
$$

 \ We will use now the following elementary inequality

$$
x^r \le \left( \frac{r}{\lambda \ e} \right)^r \cdot e^{(\lambda, \ x)}, \hspace{3mm} r, \ \lambda, \ x > 0. \eqno(7.13)
$$

 \ As a consequence: let $  \xi $ be non - zero  $ d  - $  dimensional mean zero random vector belonging to the space $  B(\phi).  $
Then

$$
{\bf E} |\xi|^r \le 2^d \  \left( \frac{r}{\lambda \ e} \right)^r  \ e^{ \phi(\lambda ||\xi||) } =
2^d \ e^{-|r|} \ r^r \ \lambda^{-r} \ e^{ \phi(\lambda ||\xi||) }, \ \lambda > 0. \eqno(7.14)
$$

 \ We find likewise the one-dimensional case: \\

\vspace{4mm}

{\bf Proposition 7.2.} Let $ \phi (\cdot) $  be arbitrary non - negative continuous function and let the centered numerical
r.v. $  \xi $  be such that $  \xi \in B(\phi): \ 0 < ||\xi|| = ||\xi||B(\phi) < \infty. $ \par

 \ Then

$$
| \vec{\xi}|_{\vec{r}} \le e^{-1} \cdot  2^{d/|r|} \cdot \prod_j r(j)^{r(j)/|r|}  \cdot  e^{- \Phi^*(r)/|r| } \cdot ||\xi||B(\phi),
 \ r = \vec{r} > 0. \eqno(7.15)
$$

 \ Note that in general case the expression $ |\xi|_r  $ does not represent the norm relative the random vector $  \vec{\xi}. $ \par
 \ But if we denote

$$
\psi_{\Phi}(\vec{r}) :=  e^{-1} \cdot 2^{d/|r|} \cdot \prod_j r(j)^{r(j)/|r|}  \cdot  e^{- \Phi^*(r)/|r| }
$$
  and define
$$
||\xi||G\psi_{\Phi} \stackrel{def}{=} \sup_{\vec{r} \ge 1} \left[  \frac{| \vec{\xi}|_{\vec{r}} }{\psi_{\Phi}(\vec{r}) }  \right],
$$
we obtain some modification of the one-dimensional Grand Lebesgue Space (GLS) norm, see (1.6). \par

 \ The statement of proposition 7.2 may be rewritten as follows.

$$
||\xi||G\psi_{\Phi}  \le ||\xi||B(\phi). \eqno(7.15a)
$$

 \ Let us state the inverse up to multiplicative constant inequality. \par

\vspace{4mm}

{\bf Theorem 7.2.} Suppose $  \phi \in S(R^d). $  Let $ {\bf E} \xi =0  $ and
$ K = ||\xi||G\psi_{\Phi} < \infty. $ Then $ \xi \in B(\phi)  $ and moreover both the norms
$ || \xi||B(\phi) $  and  $ ||\xi||G\psi_{\Phi} $ are equivalent:

$$
|| \xi||B(\phi) \le C_3(\phi) \ ||\xi||G\psi_{\Phi}, \   C_3(\phi)\in (0,\infty).  \eqno(7.16)
$$

 \ The proof of this theorem is at the same as one in the theorem  7.1 and may be omitted. \par

 \vspace{3mm}

\section{Exponential bounds for the sums of random vectors.}

\vspace{3mm}

 \ Statement of problem: given independent centered   $  d \ - $ dimensional random
vectors $ \xi_i, \ i = 1,2,\ldots;  $ put as above

$$
S(n) = n^{-1/2} \sum_{i=1}^n \xi_i,
$$
the classical norming for the sum of independent mean zero i.d. random vectors. \par

 \ It is required to deduce the exponential estimate for tails $  U(S(n), \ \vec{x}), $ (non-uniform estimates),
as well as for uniform tails   $  \sup_n U(S(n), \ \vec{x}), $  for sufficiently greatest values $ \vec{x}.  $ \par

\vspace{3mm}

 \ We claim to obtain an upper as well as lower exponential decreasing bounds for mentioned probabilities. \par

\vspace{3mm}

\ Of course, these estimates may be used in statistics and in the Monte - Carlo method, see \cite{Frolov1},
\cite{Grigorjeva1}, namely, in the so-called "method of depending trials". \par

\vspace{3mm}

 \ Related previous works: \cite{Prokhorov1}, \cite{Kurbanmuradov1}, \cite{Ostrovsky1}, chapters 2,3;
 \cite{Ostrovsky5}. \par

\vspace{3mm}

{\bf A. Upper estimates.} \par

\vspace{3mm}

 \ {\bf Definition 8.1.}  The function $  \phi: R^d \to R $ from the set $  Y = Y(R^d)  $ belongs by definition to the
class $  \Lambda_2,\ \phi(\cdot) \in \Lambda_2, $  iff for all positive numbers $  a,b > 0  $ and for all the vectors $  \lambda \in R^d  $

$$
\phi(a \cdot \lambda) + \phi(b \cdot \lambda) \le \phi(  \sqrt{ a^2 + b^2} \cdot \lambda ). \eqno(8.1)
$$

\vspace{3mm}

 \ This condition can be rewritten as follows: $  \forall c = \const > 0  $ and for all the vectors
 $  \lambda \in R^d  \ \Rightarrow $

 $$
\phi(\lambda) + \phi(c \cdot \lambda) \le \phi(  \sqrt{ c^2 + 1} \cdot \lambda ) \eqno(8.1a)
$$
or  also  as follows: $  \forall \theta \in [0, \pi/2] \ $ and for all the vectors
$  \lambda \in R^d  \ \ \Rightarrow  $

 $$
\phi(\lambda \ \cos (\theta) ) + \phi(  \lambda \ \sin(\theta) ) \le \phi( \lambda ). \eqno(8.1b)
$$

 \ The condition (8.1) is trivially satisfied if for example

$$
\phi(\lambda) =  (Q\lambda,\lambda),
$$
 where  $  Q  $ is symmetrical positive definite  matrix of a size $ d \times d. $ \par

 \vspace{3mm}

 \ {\bf Lemma 8.1.} Assume that the  function $  \phi(\cdot) $ has a form

$$
\phi(\lambda)   =  \nu((Q\lambda,\lambda)), \eqno(8.2)
$$
where  $  Q  $ is symmetrical positive definite  matrix of a size $ d \times d, $
 $  \nu: R^1 \to R^1  $ is continuous numerical
non - negative convex function such that $  \forall z \ge 0 \ \nu(z) = 0 \Leftrightarrow \ z = 0. $ \par

\ We propose that  $ \phi(\cdot)  \in \Lambda_2. $ \par

 \vspace{3mm}

  \ {\bf  An example: }

 $$
 \phi(\lambda) = \nu( \ |\lambda|^2 \ ), \eqno(8.3)
 $$
where the properties of the function $  \nu(\cdot)  $  was described  before. \par

 \ The functions of the form (8.3) are named spherical, or radial. \par

\vspace{3mm}

 \ {\bf Proof.}  We will use the following elementary for such a convex functions $  \{  \nu(\cdot) \}  $

$$
\nu(x) + \nu(y) \le \nu(x + y), \ x,y \ge 0,
$$
see \cite{Krasnosel'skii1}, chapter 1, sections 1,2.  Therefore

$$
\nu(a^2 (Q \lambda,\lambda)) +  \nu(b^2 (Q \lambda,\lambda)) \le  \nu( (a^2 + b^2) \cdot (Q \lambda,\lambda)),
$$
which is equivalent to the inequality (8.1). \par

\vspace{3mm}

 \ {\bf Theorem 8.1.} {\it Suppose the function $ \phi = \phi(\lambda)  $ belongs to the set $  \Lambda_2.  $ If
 $  \{\xi_i\}, \ i = 1,2,\ldots,n   $ are independent random vectors belonging to the space $  B(\phi)  $ with
 finite norms in this space $ ||\xi_i||B(\phi) =  ||\xi_i||.  $ We propose }

$$
||S(n)||B(\phi) \le n^{-1/2} \left[ \sum_{i=1}^n  ||\xi_i||^2 \right]^{1/2} \stackrel{def}{=} \sigma(n). \eqno (8.4)
$$

\vspace{3mm}

 \ {\bf Proof} follows immediately from the estimate

$$
||\xi + \eta|| \le ||\xi|| \odot_{\phi} ||\eta|| \le \sqrt{ ||\xi||^2 + ||\eta||^2}, \eqno(8.5)
$$
which take place for independent random vectors $  \xi, \ \eta   $  from the space $  B(\phi), $
(Pythagoras inequality), and from the homogeneity of the  norm   $ ||\xi||B(\phi).  $\par

\vspace{3mm}

 \ It follows immediately from the last inequality (6.5) and theorem 6.1 \\

\vspace{3mm}

 \ {\bf Corollary 8.1.} {\it We  affirm under  conditions of theorem 8.1}

$$
 U(S(n), \ \vec{x})  \le \exp \left( - \phi^*(\vec{x}/\sigma(n)) \right), \ \vec{x} \ge 0. \eqno(8.6)
$$

 \ {\it If in addition } $   \sigma := \sup_{n} \sigma(n) < \infty, $ {\it then of course}

$$
\sup_n U(S(n), \ \vec{x})  \le \exp \left( - \phi^*(\vec{x}/\sigma) \right). \eqno(8.6a)
$$

\vspace{3mm}

 \ Another approach. Let the r.v. $ \vec{\xi} = \xi  $  be from some space $  B(\phi), \ \phi \in Y(V). $
For instance, the function $  \phi(\cdot) $ may be the natural function for the r.v. $  \xi: \ \phi(\lambda) =
\phi_{\xi}(\lambda). $ \par

  \ Introduce a new functions

$$
\phi_n(\lambda) \stackrel{def}{=} n \ \phi(\lambda/\sqrt{n}), \hspace{5mm}
\overline{\phi}(\lambda) \stackrel{def}{=} \sup_n  [ n \ \phi(\lambda/\sqrt{n})].
 \eqno(8.7)
$$

 \ The last function $ \overline{\phi}(\lambda), \ \lambda \in \Lambda  $ is finite as long as there exists a limit

 $$
 \lim_{n \to \infty} n  \phi(\lambda/\sqrt{n}) = 0.5 \ (\phi^{''}(0) \lambda, \lambda) < \infty.
 $$

 \ Let us explain the sense of this function. Suppose here the random vectors $ \{  \xi_i \}, \ i = 1,2,\ldots,n   $
 are independent identically distributed copies of $  \xi,  $  then

$$
{\bf E} e^{\lambda S(n)} \le e^{ \phi_n(\lambda) }, \eqno(8.8)
$$
and  correspondingly

$$
\sup_n {\bf E} e^{\lambda S(n)} \le e^{ \overline{\phi} (\lambda) }. \eqno(8.8a)
$$

\vspace{3mm}

 {\bf  Proposition 8.1.} We conclude alike the theorem 8.1 applying theorem 6.1 under formulated in this
 pilcrow conditions

$$
 U(S(n), \ \vec{x})  \le \exp \left( - \phi_n^*(\vec{x}) \right), \ \vec{x} \ge 0. \eqno(8.9)
$$

$$
\sup_n U(S(n), \ \vec{x})  \le \exp \left( - \overline{\phi}^*(\vec{x} \right), \ \vec{x} \ge 0. \eqno(8.9a)
$$

\vspace{3mm}

{\bf B. \ Lower bounds.} \par

\vspace{3mm}

 \ Lower estimates for the considered here probabilities are simple; we follow \cite{Ostrovsky1}, p. 50 - 53,
where it was considered the one-dimensional case $  d = 1. $\par

 \ Namely, let here $  \xi, \ \{ \xi_i \}, \ i = 1,2,\ldots $ be a sequence of centered identical distributed independent
  random vectors having exponential decreasing tails, i.e. belonging  to some space $  B(\phi),  $ where $  \phi(\cdot) \in S(R^d). $
 Then we have on the  one hand

 $$
\sup_n U(S(n), \ \vec{x})  \ge U( \vec{\xi}, \vec{x}). \eqno(8.10)
 $$

  \ On the  one hand,

$$
\sup_n U(S(n), \ \vec{x})  \ge \lim_{n \to \infty} U(S(n), \ \vec{x}). \eqno(8.11)
$$

 \ Let us denote  $ Q = \Var(\xi)   $ and by $  \gamma_Q(\cdot)  $ the Gaussian distribution (measure) on the
 space $   R^d $ with mean zero and variance $  Q. $ Then by virtue of multivariate CLT

$$
\lim_{n \to \infty} U(S(n), \ \vec{x}) = \int_{ \vec{y} \ge \vec{x} } \gamma_Q(d y). \eqno(8.12)
$$

 \ The right - hand side of the relation (8.12) may be estimated as follows:

$$
\int_{ \vec{y} \ge \vec{x} } \gamma_Q(d y) \ge \exp \left( - C(Q) \ |x|^2  \right), \  |x|^2 = (x,x) \ge 1.
$$

 \ We established in fact the following lower estimate. \par

 \vspace{3mm}

 \ {\bf Proposition 8.2.}

 $$
\sup_n U(S(n), \ \vec{x})  \ge \max \left( U( \vec{\xi}, \vec{x}), \ \exp \left( - C(Q) \ |x|^2  \right) \ \right), \  |x| \ge 1.
 \eqno(8.13)
 $$

 \vspace{3mm}

{\bf Example 8.1.} Let as before $  \xi, \ \{ \xi_i \}, \ i = 1,2,\ldots $ be a sequence of centered identical
distributed independent random vectors such that

$$
 U( \vec{\xi}, \vec{x}) \le  \ e^{ - |x|^p }, \ |x| \ge 1,
$$
where $  p = \const \ge 1. $ Then

$$
\sup_n U(S(n), \ \vec{x}) \le \exp  \left( - \ C(d,p) \ |\vec{x}|^{\min(p,2)} \ \right),   \   |x| \ge 1.
$$

\vspace{3mm}

{\bf Remark 8.1.}  \ Note  that the last estimate can not be improvable still in the one-dimensional case, see
\cite{Kozachenko1},  \cite{Ostrovsky1}, p. 50-53. \par

\vspace{3mm}

{\bf Remark 8.2.} \ Note finally that the last estimate can not be obtained by means of Rosenthal's moment inequality,
although the $  B(\phi)  $  and moment (GLS) norm are (linear) equivalent.\par

 \ This is due to the fact that the famous Rosenthal's "constants"  $ R(p) $ tending to  infinity  as $  p \to \infty:  $

$$
R(p) \asymp \frac{p}{\log p}, \ p \in [2, \infty).
$$

 \ The detail explanation of this phenomenon may be found in \cite{Ostrovsky1}, chapter 2, sections 2.1 and 2.3.\par

\vspace{3mm}

\section{Relation with multivariate Orlicz spaces.}

\vspace{3mm}

 \ Let now the function $  \phi = \phi(\lambda) $ be from the set $  Y(V), $ see definition 2.5. Introduce
the following Young-Orlicz function

$$
N_{\phi}(u) = N_{\phi}(\vec{u}) = \exp (\phi^*(u)) - \exp (\phi^*(0)), \eqno(9.1)
$$
so that $  N(u) $ is even function and $  N(0) = 0.  $ \par

 \ An Orlicz space on the source probability triple generated by the function $ N_{\phi}(\vec{u}) $ will be denoted as
ordinary by $  L(N_{\phi})  $ and the correspondent norm by $  ||\xi||L(N_{\phi}).  $ \par

\vspace{3mm}

 \ {\bf  Theorem 9.1. } Both the norms $  ||\xi||L(N_{\phi})  $ and $   ||\xi||B(\phi) $ are equivalent on the subspace of the
 centered random vectors:

$$
\exists \ C_{5,6} = C_{5,6}(d, \phi) \in (0,\infty),  \ 0 < C_5 < C_6 < \infty
$$
such that for all the $ d \ - $ dimensional random vectors $ \xi: \ {\bf E} \xi = 0  $

$$
C_5 ||\xi||L(N_{\phi})  \le  ||\xi||B(\phi) \le   C_6 ||\xi||L(N_{\phi}). \eqno(9.2)
$$

\vspace{3mm}

 \ {\bf Proof} is at the same as in the one - dimensional case, see \cite{Kozachenko1}. Suppose first of all
 $  {\bf E} \xi = 0  $ and $   ||\xi||L(N_{\phi}) = 1/2.  $ Then

$$
{\bf E} N_{\phi} (\xi) \le 1.
$$

  \ The Tchebychev's  inequality allows us to estimate the tail behavior of the r.v. $  \xi  $ as follows:

$$
U(\vec{\xi}, \vec{x}) \le \frac{1}{N_{\phi}(\vec{x})} \le \exp \left( - \phi^*(c \ \vec{x})  \right), \ |x| \ge 1.
$$
  The case $ |x| \le 1 $ is  obvious.  The right - hand side of bilateral inequality (9.2) follows immediately from
theorem 6.1. \par

 \ Let us prove now the inverse estimate. Suppose for instance

$$
{\bf E} e^{ (\lambda,\xi) } \le e^{\phi(\lambda)}, \ \lambda \in R^d_+.
$$
 Then $ {\bf E} \vec{\xi} = 0 $ and

$$
U(\vec{\xi}, \vec{x}) \le e^{ - \phi^*(\vec{x}) }, \ \vec{x} > 0.
$$

 \ We have analogously to the preprint \cite{Ostrovsky5}

$$
{\bf E} N( 0.5 \ \vec{\xi} ) \le C_1 +  C_2 \int_{R^d_+} \exp \left( \phi^*(0.5 \ x) - \phi^*(x) \right) \ dx < \infty,
$$
therefore $  \vec{\xi} \in L(N_{\phi}).  $ \par

\vspace{3mm}

\section{Concluding remarks.}

\vspace{3mm}

{\bf A.  Another approach to the problem of description of natural function.} \par

\vspace{3mm}

 \ Recall the statement of problem, see fourth section. Given a function $  \phi(\cdot) $  from the set $ Y = Y(R^d). $ Question: under what
additional conditions there exists a mean zero random   vector $  \xi  $ which may be defined on the appropriate
probability space such that

$$
\forall \lambda \in R^d \ \Rightarrow {\bf E} e^{ (\lambda, \ \xi) } = e^{\phi(\lambda)}. \eqno(10.1)
$$

 \ Obviously, the function $  \lambda \to e^{\phi(\lambda)} $ may be analytically  continued  on the whole complex space $  {\bf C^d } $
 such that

$$
e^{\phi( i \cdot \vec{t})} = {\bf E} e^{i ( \ \vec{t}, \ \vec{\xi}) }, \ i^2 = -1, \ \vec{t} \in R^d, \eqno(10.2)
$$
so that $ \exp \{\phi( i \cdot \vec{t} )  \} $ is a characteristical function of some $ d \ - $ dimensional random vector. \par

\vspace{4mm}

{\bf B. Absolutely even functions.} \par

\vspace{3mm}

 \ The  function $ f: R^d \to R   $  is said to be absolutely even, in $  \forall \vec{\epsilon}, \ \vec{x} $

$$
f(\vec{\epsilon} \otimes \vec{x}) \ = f( \vec{x}). \eqno(10.3)
$$

 \ In other words, it is even function separately arbitrary variable by fixed others. \par

 \ For instance, the function of a form

$$
f(x,y) = g(x,y) + g(x,-y) + g(-x,y) + g(-x, -y)
$$
as well as

$$
f_1(x,y) = \max[ g(x,y), g(x,-y), g(-x,y), g(-x, -y) ]
$$
are absolutely even; here $ d = 2 $ and $  g: R^2 \to R.  $ \par

 \ It is interest to note by our opinion that the natural function for any random vector is absolutely even.\par

\vspace{4mm}

{\bf C. Conjugate and associate spaces.} \par

\vspace{3mm}

 \ Since the introduced spaces $  B( \vec{\phi})  $ coincide with Orlicz spaces, see the $  9^{th} $ section, the
study of its conjugate and associate spaces may be provided likewise as in monographs  \cite{Rao1}, \cite{Rao2}. \par
 See also \cite{Ostrovsky5}. \par

\end{document}